\documentclass{article}%
\usepackage{graphicx}
\usepackage{amsmath}
\usepackage{amsfonts}
\usepackage{amssymb}%
\newtheorem{theorem}{Theorem}

\newtheorem{corollary}[theorem]{Corollary}

\begin{document}

\title{A remark on absolutely summing multilinear mappings}
\author{Daniel M. Pellegrino\thanks{Partially supported by Instituto do Milenio}\\DME-UFCG\\Caixa Postal 10044\\Campina Grande, PB, 58109-970 Brazil}
\date{e-mail address: dmp@dme.ufcg.edu.br }
\maketitle

\begin{abstract}
In this note we obtain new coincidence theorems for absolutely summing
multilinear mappings between Banach spaces. We also prove that our results, in
general, can not be improved.

\end{abstract}



\section{Introduction}

The theory of ideals of multilinear mappings was outlined by
Pietsch \cite{Pietsch} and since then the generalization of
particular ideals of operators to multilinear mappings has being
investigated by several authors (see
\cite{BBJ},\cite{Cilia},\cite{Dimant},\cite{Matos3},\cite{P1},\cite{P2}).
The analogy between the linear and multilinear cases is sometimes
non trivial and several related questions have been investigated.
For example, it is well known that every scalar valued continuous
linear operator between Banach spaces is absolutely $p$-summing.
It is also easy to check that this result, in general, is no
longer valid for absolutely $(p;p,...,p)$-summing multilinear
mappings. However, an unpublished result, credited to A. Defant
and J. Voigt, states that every continuous multilinear form is
absolutely $(1;1,...,1)$-summing (see \cite{Matos}). In this note
we obtain new results of coincidence for absolutely summing
multilinear mappings and show that, in general, they can not be
improved.

Throughout this paper $p$ is a real number not smaller than $1$ and $E$,
$E_{1}$,..., $E_{n}$ and $F$ are Banach spaces. The scalar field $\mathbb{K}$
can be either $\mathbb{R}$ or $\mathbb{C}.$ The linear space of all sequences
$(x_{j})_{j=1}^{\infty}$ in $E$ such that $\Vert(x_{j})_{j=1}^{\infty}%
\Vert_{p}=(\sum_{j=1}^{\infty}\Vert x_{j}\Vert^{p})^{\frac{1}{p}}<\infty$ will
be represented by $l_{p}(E).$ We will also denote by $l_{p}^{w}(E)$ the linear
space composed by the sequences $(x_{j})_{j=1}^{\infty}$ in $E$ such that
$(\varphi(x_{j}))_{j=1}^{\infty}\in l_{p}$ for every bounded linear functional
$\varphi:E\rightarrow\mathbb{K}.$ We define $\Vert.\Vert_{w,p}$ in $l_{p}%
^{w}(E)$ by $\Vert(x_{j})_{j=1}^{\infty}\Vert_{w,p}=\sup_{\varphi\in
B_{E^{\prime}}}\Vert(\varphi(x_{j}))_{j=1}^{\infty}\Vert_{p}.$

The Banach space of all bounded $n$-linear mappings from $E_{1}\times...\times
E_{n}$ into $F$ endowed with the $\sup$ norm will be represented by
$\mathcal{L}(E_{1},...,E_{n}$;$F)$ and the Banach space of all continuous
$n$-homogeneous polynomials $P$ from $E$ into $F$ with the $\sup$ norm is
denoted by $\mathcal{P}(^{n}E;F).$

The definition of absolutely summing polynomials and multilinear mappings we
will work with appears in \cite{Alencar} and is a natural generalization of
the linear case:

$T\in\mathcal{L}(E_{1},...,E_{n};F)$ is absolutely $(p;q_{1},...,q_{n}%
)$-summing if $(T(x_{j}^{(1)},...,x_{j}^{(n)}))_{j=1}^{\infty}\in l_{p}(F)$
for all $(x_{j}^{(s)})_{j=1}^{\infty}\in l_{q_{s}}^{w}(E_{s})$, $s=1,...,n.$
We will write $\mathcal{L}_{as(p;q_{1},...,q_{n})}(E_{1},...,E_{n};F)$ to
denote the space of all absolutely $(p;q_{1},...,q_{n})$-summing multilinear
mappings from $E_{1}\times...\times E_{n}$ into $F.$ As in the linear case, it
is known that $T$ $\in\mathcal{L}_{as(p;q_{1},...,q_{n})}(E_{1},...,E_{n};F)$
if and only if there exists $C>0$ such that
\[
(\sum_{j=1}^{k}\Vert T(x_{j}^{(1)},...,x_{j}^{(n)})\Vert^{p})^{\frac{1}{p}%
}\leq C\Vert(x_{j}^{(1)})_{j=1}^{k}\Vert_{w,q_{1}}...\Vert(x_{j}^{(n)}%
)_{j=1}^{k}\Vert_{w,q_{n}}%
\]
for each natural $k$ and $x_{j}^{(l)}\in E_{l},$ $l=1,...,n$ and $j=1,...,k$
(see Matos \cite{Matos})$.$

The infimum of the $C>0$ for which the last inequality holds defines a norm
$\Vert.\Vert_{as(p;q_{1},...,q_{n})}$ for $\mathcal{L}_{as(p;q_{1},...,q_{n}%
)}(E_{1},...,E_{n};F),$ and under this norm this space is complete.

\section{Results}

In \cite{Botelho}, Botelho shows that every continuous bilinear form defined
in $\mathcal{L}_{\infty}$-spaces (for the definition of $\mathcal{L}_{\infty}%
$-spaces we mention \cite{Lindenstrauss}) is absolutely $(1;2,2)$-summing
($2$-dominated)$.$ In the same paper it is also proved that we can not expect
another similar coincidence theorem for $p$-dominated $n$-linear mappings,
with $n>2$. Recently, using a generalized Grothendieck's inequality,
P\'{e}rez-Garc\'{\i}a \cite{Perez} obtained the following result of coincidence:

\begin{theorem}
(P\'{e}rez-Garc\'{\i}a \cite{Perez})\label{tp} If $E_{1},...,E_{n}$ are
$\mathcal{L}_{\infty}$-spaces, then every continuous $n$-linear ($n\geq2$)
mapping $T:E_{1}\times...\times E_{n}\rightarrow\mathbb{K}$ is absolutely
$(1;2...,2)$-summing.
\end{theorem}

Our first goal is to obtain coincidence results for the case in
which some of the $E_{j}$ are $\mathcal{L}_{\infty}$-spaces  and
some are arbitrary Banach spaces. The next result generalizes a
theorem of C.A. Soares \cite{CA} and is crucial for our purposes.

\begin{theorem}
\label{t1}Let $A\in\mathcal{L}(E_{1},...,E_{n};F).$ Suppose that there exists
$K>0$ so that for any $x_{1}\in E_{1},....,x_{r}\in E_{r},$ the $s$-linear
($s=n-r$) mapping $A_{x_{1}....x_{r}}(x_{r+1},...,x_{n})=A(x_{1},...,x_{n})$
is absolutely $(1;q_{1},...,q_{s})$-summing and besides $\left\|
A_{x_{1}....x_{r}}\right\|  _{as(1;q_{1},...,q_{s})}\leq K\left\|  A\right\|
\left\|  x_{1}\right\|  ...\left\|  x_{r}\right\|  $. Then we can conclude
that $A$ is absolutely $(1;1,...,1,q_{1},...,q_{s})$-summing.
\end{theorem}

Proof: For $x_{1}^{(1)},...,x_{1}^{(m)}\in E_{1},....,x_{n}^{(1)}%
,...,x_{n}^{(m)}\in E_{n}$, let us consider $\varphi_{j}\in B_{F^{\prime}}$
such that $\left\Vert A(x_{1}^{(j)},...,x_{n}^{(j)})\right\Vert =\varphi
_{j}(A(x_{1}^{(j)},...,x_{n}^{(j)}))\text{ for every }j=1,...,m.$ Thus,
defining by $r_{j}(t)$ the Rademacher functions on $[0,1]$ and denoting by
$\lambda$ the Lebesgue measure in $I=[0,1]^{r},$ we have
\begin{align*}
&  \int\nolimits_{I}\sum\limits_{j=1}^{m}\left(  \prod_{l=1}^{r}r_{j}%
(t_{l})\right)  \varphi_{j}A(\sum\limits_{j_{1}=1}^{m}r_{j_{1}}(t_{1}%
)x_{1}^{(j_{1})},...,\sum\limits_{j_{r}=1}^{m}r_{j_{r}}(t_{r})x_{r}^{(j_{r}%
)},x_{r+1}^{(j)},...,x_{n}^{(j)})d\lambda\\
&  =\sum\limits_{j,j_{1},...j_{r}=1}^{m}\varphi_{j}A(x_{1}^{(j_{1})}%
,...,x_{r}^{(j_{r})},x_{r+1}^{(j)},...,x_{n}^{(j)})\int\limits_{0}^{1}%
r_{j}(t_{1})r_{j_{1}}(t_{1})dt_{1}...\int\limits_{0}^{1}r_{j}(t_{r})r_{j_{r}%
}(t_{r})dt_{r}%
\end{align*}%
\[
=\sum\limits_{j=1}^{m}\varphi_{j}A(x_{1}^{(j)},...,x_{n}^{(j)})=\sum
\limits_{j=1}^{m}\left\Vert A(x_{1}^{(j)},...,x_{n}^{(j)})\right\Vert
=(\ast).
\]
So, for each $l=1,...,r,$ assuming $z_{l}=\sum\limits_{j=1}^{m}r_{j}%
(t_{l})x_{l}^{(j)}$ we obtain
\begin{align*}
(\ast)  &  =\int\nolimits_{I}\sum\limits_{j=1}^{m}\left(  \prod_{l=1}^{r}%
r_{j}(t_{l})\right)  \varphi_{j}A(\sum\limits_{j_{1}=1}^{m}r_{j_{1}}%
(t_{1})x_{1}^{(j_{1})},...,\sum\limits_{j_{r}=1}^{m}r_{j_{r}}(t_{r}%
)x_{r}^{(j_{r})},x_{r+1}^{(j)},...,x_{n}^{(j)})d\lambda\\
&  \leq\int\nolimits_{I}\left\vert \sum\limits_{j=1}^{m}\left(  \prod
_{l=1}^{r}r_{j}(t_{l})\right)  \varphi_{j}A(\sum\limits_{j_{1}=1}^{m}r_{j_{1}%
}(t_{1})x_{1}^{(j_{1})},...,\sum\limits_{j_{r}=1}^{m}r_{j_{r}}(t_{r}%
)x_{r}^{(j_{r})},x_{r+1}^{(j)},...,x_{n}^{(j)})\right\vert d\lambda\\
&  \leq\int\nolimits_{I}\sum\limits_{j=1}^{m}\left\Vert A(\sum\limits_{j_{1}%
=1}^{m}r_{j_{1}}(t_{1})x_{1}^{(j_{1})},...,\sum\limits_{j_{r}=1}^{m}r_{j_{r}%
}(t_{r})x_{r}^{(j_{r})},x_{r+1}^{(j)},...,x_{n}^{(j)})\right\Vert d\lambda\\
&  \leq\sup_{t_{l}\in\lbrack0,1],l=1,...,r}\sum\limits_{j=1}^{m}\left\Vert
A(\sum\limits_{j_{1}=1}^{m}r_{j_{1}}(t_{1})x_{1}^{(j_{1})},...,\sum
\limits_{j_{r}=1}^{m}r_{j_{r}}(t_{r})x_{r}^{(j_{r})},x_{r+1}^{(j)}%
,...,x_{n}^{(j)})\right\Vert \\
&  \leq\sup_{t_{l}\in\lbrack0,1],l=1,...,r}\left\Vert A_{z_{1}...z_{r}%
}\right\Vert _{as(1;q_{1},...,q_{s})}\left\Vert (x_{r+1}^{(j)})_{j=1}%
^{m}\right\Vert _{w,q_{1}}...\left\Vert (x_{n}^{(j)})_{j=1}^{m}\right\Vert
_{w,q_{s}}\\
&  =K\left\Vert A\right\Vert \left(  \prod_{l=1}^{r}\sup_{t\in\lbrack
0,1]}\left\Vert \sum\limits_{j=1}^{m}r_{j}(t)x_{l}^{(j)}\right\Vert \right)
\left(  \prod_{l=1}^{s}\left\Vert (x_{l}^{(j)})_{j=1}^{m}\right\Vert
_{w,q_{l}}\right) \\
&  \leq K\left\Vert A\right\Vert \left(  \prod_{l=1}^{r}\left\Vert
(x_{l}^{(j)})_{j=1}^{m}\right\Vert _{w,1}\right)  \left(  \prod_{l=1}%
^{s}\left\Vert (x_{l}^{(j)})_{j=1}^{m}\right\Vert _{w,q_{l}}\right)  .\text{
Q.E.D.}%
\end{align*}
We have the following straightforward consequence:

\begin{corollary}
\label{ssss}If $\mathcal{L}(E_{1},...,E_{m};F)=\mathcal{L}_{as(1;q_{1}%
,...,q_{m})}(E_{1},...,E_{m};F)$ then, for any Banach spaces $E_{m+1}%
,...,E_{n},$ we have $\mathcal{L}(E_{1},...,E_{n};F)=\mathcal{L}%
_{as(1;q_{1},...,q_{m},1,...,1)}(E_{1},...,E_{n};F).$
\end{corollary}

Note that an particular case of this result is the aforementioned coincidence
result of Defant and Voigt. Another outcome of Theorem \ref{t1} is the
following corollary, which proof is an immediate consequence of Corollary
\ref{ssss} and Theorem \ref{tp}:

\begin{corollary}
\label{c1}If $E_{1}$,..., $E_{s}$ are $\mathcal{L}_{\infty}$-spaces then, for
any Banach spaces $E_{s+1},...,E_{n},$ we have $\mathcal{L}(E_{1}%
,...,E_{n};\mathbb{K})=\mathcal{L}_{as(1;q_{1},...,q_{n})}(E_{1}%
,...,E_{n};\mathbb{K})$, where $q_{1}=...=q_{s}=2$ e $q_{s+1}=....=q_{n}=1.$
\end{corollary}

The next straightforward corollaries, show that we still have interesting
coincidence results if the range of our mappings has finite cotype.

\begin{corollary}
\label{ff}If $\cot F=q<\infty$ \ and $\mathcal{L}(E_{1},...,E_{s}%
;\mathbb{K})=\mathcal{L}_{as(1;q_{1},....,q_{s})}(E_{1},...,E_{s}%
;\mathbb{K}),$ then, for any $E_{s+1},...,E_{n},$ we have $\mathcal{L}%
(E_{1},...,E_{n};F)=\mathcal{L}_{as(q;q_{1},....,q_{s},1,....,1)}%
(E_{1},...,E_{n};F).$
\end{corollary}

\begin{corollary}
\label{g}If $\cot F=q<\infty$ and $E_{1}$,..., $E_{s}$ are $\mathcal{L}%
_{\infty}$-spaces, then, regardless of the Banach spaces $E_{s+1},...,E_{n},$
we have $\mathcal{L}(E_{1},...,E_{n};F)=\mathcal{L}_{as(q;q_{1},...,q_{n}%
)}(E_{1},...,E_{n};F),$where $q_{1}=...=q_{s}=2$ and $q_{s+1}=...=q_{n}=1.$
\end{corollary}

It is obvious that Corollary \ref{c1} is still true if we replace $\mathbb{K}$
by any finite dimensional Banach space. A natural question is whether
Corollary \ref{c1} can be improved for some infinite dimensional Banach space
in the place of $\mathbb{K}.$ Precisely, the question is:

\begin{itemize}
\item If $E_{1},...,E_{k}$ are infinite dimensional $\mathcal{L}_{\infty}%
$-spaces, is there some infinite dimensional $F$ so that $\mathcal{L}%
(E_{1},...,E_{k},...,E_{n};F)=\mathcal{L}_{as(1;q_{1},....,q_{n})}%
(E_{1},...,E_{k},...,E_{n};F),$ where $q_{1}=...=q_{k}=2$ and $q_{k+1}%
=....=q_{n}=1,$ regardless of the Banach spaces $E_{k+1},...,E_{n}$?
\end{itemize}

The answer to this question, surprisingly, is no. The proof is obtained as an
application of \cite[Theorem 8]{Pellegrino}. Another relevant question is
whether Corollary \ref{g} can be improved to $p<q$, i.e.,

\begin{itemize}
\item If $\cot F=q<\infty$ and $E_{1}$,..., $E_{k}$ are infinite dimensional
$\mathcal{L}_{\infty}$-spaces, is there some $p<q$ for which, regardless of
the Banach spaces $E_{k+1},...,E_{n},$ $\mathcal{L}(E_{1},...,E_{k}%
,...,E_{n};F)=\mathcal{L}_{as(p;q_{1},....,q_{n})}(E_{1},...,E_{k}%
,...,E_{n};F),$ where $q_{1}=...=q_{k}=2$ , $q_{k+1}=....=q_{n}=1$?
\end{itemize}

Again, using \cite[Theorem 8]{Pellegrino} one can obtain a negative answer to
this question.

This note forms a part of the author's thesis written under supervision of
M.C. Matos, and he wishes to thank him for the suggestions.


\begin{thebibliography}{99}                                                                                               %


\bibitem {Alencar}R. Alencar and M. C. Matos, Some classes of multilinear
mappings between Banach spaces, Publ. Dep. An\'{a}lisis Mat. Univ. Complut.
Madrid 12 (1989).

\bibitem {Botelho}G. Botelho, Cotype and absolutely summing multilinear
mappings and homogeneous polynomials, Proceedings of The Royal Irish Academy
Vol 97A, No 2 (1997), 145-153.


\bibitem {BBJ}G. Botelho, H. A. Braunss and H. Junek, Almost $p$-summing
polynomials and multilinear mappings, Arch. Math. (Basel)
\textbf{76} (2001), 109-118.

\bibitem {Cilia}R. Cilia, M. D'Anna and J. Guti\'{e}rrez, Polynomial
characterization of $\mathcal{L}_{\infty}$-spaces, J. Math. Anal.
Appl. \textbf{275} (2002) 900-912.

\bibitem {Dimant}V. Dimant, Strongly $p$-summing multilinear operators, J.
Math. Anal. Appl. \textbf{278} (2003), 182-193.

\bibitem {Lindenstrauss}J. Lindenstrauss and A. Pe\l czy\'{n}ski,\ Absolutely
summing operators in $\mathcal{L}_{p}$ spaces and their
applications, Studia Mathematica \textbf{29} (1968), 275-326.

\bibitem {Matos}M. C. Matos, Absolutely summing holomorphic mappings, An.
Acad. bras. Ci., \textbf{68} (1996), 1-13.

\bibitem {Matos3}M. C. Matos, Fully absolutely summing and Hilbert Schmidt
multilinear mappings, Collect. Mat. \textbf{54} (2003), 111-136.

\bibitem {Pellegrino}D. M. Pellegrino, Cotype and absolutely summing
homogeneous polynomials in $L_{p}$ spaces, Studia Math.
\textbf{157} (2003), 121-131.

\bibitem {P1}D. M. Pellegrino, Almost summing multilinear mappings, to appear
in Arch. Math (Basel).

\bibitem {P2}D. M. Pellegrino, Strongly almost summing holomorphic mappings,
to appear in J. Math. Anal. Appl.

\bibitem {Perez}D. P\'{e}rez-Garc\'{\i}a, Operadores multilineales
absolutamente sumantes, Dissertation, Universidad Complutense de Madrid (2002).

\bibitem {Pietsch}A. Pietsch, Ideals of multilinear functionals (designs of a
theory), Proceedings of the Second International Conference on Operator
Algebras, Ideals and their Applications in Theoretical Physics, 185-199.
Leipzig. Teubner-Texte (1983).

\bibitem {CA}C.A. Soares, Aplica\c{c}\~{o}es multilineares e polin\^{o}mios
misto somantes, Thesis, Unicamp (1998).
\end{thebibliography}
\end{document}